\documentclass{article}
\usepackage{graphicx} 
\usepackage{tikz}
\usepackage{amsmath}
\usepackage{amssymb}
\usepackage{float}
\usepackage{wrapfig}
\usepackage{amsthm,mathtools}
\usepackage{setspace}
\usepackage[numbers]{natbib}
\usepackage{mathrsfs}
\usepackage{multicol}
\usepackage{hyperref}

\newtheorem{definition}{\hspace{2em}Definition}
\newtheorem{theorem}{\hspace{2em}Theorem}

\newtheorem{conjecture}{\hspace{2em}Conjecture}

\newtheorem{lemma}{\hspace{2em}Lemma}
\newtheorem{corollary}{\hspace{2em}Corollary}
\newtheorem{remark}{\hspace{2em}Remark}
\newcommand{\Hh}{\mathcal H}\newcommand{\Ff}{\mathcal F}\newcommand{\Gg}{\mathcal G}

\newcommand{\lstar}{\log_2^*}

\title{A Log-Star Comparison Between Expectation Threshold and Fractional Expectation Threshold}
\author{Xuan Fang\footnote{26110180010@m.fudan.edu.cn} \quad and \quad Tianyu Wang\footnote{wangtianyu@fudan.edu.cn}}
\date{Revised on Nov. 6, first version released on Sep. 25, 2026} 

\begin{document}

\maketitle

\begin{abstract}
    We proved a log-star comparison between the expectation threshold $q(\mathcal F)$ and the fractional
    expectation threshold $q_f(\mathcal F)$ for any nontrivial increasing family $\mathcal F$ on a finite ground set $V$ of size $|V|=N$. Specifically, we show that
    $$q_f(\Ff)\le64\log_2^*(N+2)\,q(\Ff),$$
    where $\log_2^* x$ is the iterated logarithm of base-2. 

    Using similiar methods, one can extend this result to the following stronger form: there exists a universal constant $C>0$ such that
    \[
    q_f(\Ff)\le Cq(\Ff)\log_2^*(l(\Ff)),\quad q_f(\Ff)\le Cq(\Ff)\bigl(1+\log_2^*(\operatorname{VC}(\min\Ff)+2)\bigr),
    \]
    where $l(\Ff):=\max\{2,\max_{S\in\min\Ff}|S|\}$ and $\operatorname{VC}(\min\Ff)$ denotes the VC-dimension of $\min\Ff$. 



    On October 7, 2026, OpenAI released a result proving the equivalence of fractional and integral thresholds, which supersedes the findings in this note. The authors are keeping this note for historical reference.

\end{abstract}

\section{Introduction}
The comparison between expectation and fractional expectation concerns how much weight one need to add when replacing an fractional cover of an increasing set family with a integral cover. 

Let $V$ be a finite set, write $N=|V|$, and let $\Ff\subseteq2^V$ be nontrivial and increasing. For $0<p<1$, call $\Ff$ \emph{$p$-small} if there is a family $\Gg\subseteq2^V$ such that every $I\in\Ff$ contains some $G\in\Gg$ and $\sum_{G\in\Gg}p^{|G|}\le1/2$. It is \emph{weakly $p$-small} if there are nonnegative weights $(\lambda_S)_{S\subseteq V}$ satisfying
\[
\sum_{S\subseteq I}\lambda_S\ge1\quad(I\in\Ff),
\qquad
\sum_{S\subseteq V}\lambda_Sp^{|S|}\le\frac12.
\]
The expectation threshold $q(\Ff)$ and fractional expectation threshold $q_f(\Ff)$ are the suprema of the parameters $p$ for which these respective conditions hold. If $p_c(\Ff)$ is the unique density at which an independent random subset of $V$ belongs to $\Ff$ with probability $1/2$, then
\[
q(\Ff)\le q_f(\Ff)\le p_c(\Ff).
\]

Khan-Kalai was the first to post this question. They conjectured in \cite{kahn2007thresholds} that there exists a universal constant $K$, such that for any ground set $V$ and any increasing family $\Ff$, let $l_0(\Ff)$ be the size of a largest minimal element of $\Ff$ and $l(\Ff)=\max\{2,l_0(\Ff)\}$, then
\begin{equation}\label{eq:khankalai}
    p_c(\Ff)\le Kq(\Ff)\log l(\Ff).
\end{equation}

Frankston, Kahn, Narayanan, and Park~\cite{frankston_kahn_narayanan_park2021} proved the fractional threshold bound $p_c(\Ff)\le Cq_f(\Ff)\log\ell$, and Park and Pham~\cite{park_pham2024_kahn_kalai} subsequently proved \eqref{eq:khankalai}. The logarithmic loss in these comparisons with $p_c$ is necessary in general, as exemplified by the appearance of perfect matchings in random graphs; see also~\cite{pham_sharp_selectors}. In particular, the second inequality gives 
\[
q_f(\Ff)\le Cq(\Ff)\log\ell,
\]
and hence an $O(\log(N+2))$ comparison. Sharpness of the comparison with $p_c$ does not establish sharpness of the comparison between $q_f$ and $q$.

More recent progress comes from selector processes and direct rounding of fractional covers. Building on the selector-process theorem of Park and Pham~\cite{park_pham2024_selectors}, Pham~\cite[Theorem~1.7]{pham_sharp_selectors} showed that a weak-$p$-smallness certificate supported on sets of size at most an integer $t\ge1$ can be rounded to an integral cover witnessing $(cp/\log(2t))$-smallness. Combining this result with the sampling reduction of Fischer and Person~\cite{fischer_person2025} gives the best previously known general bound expressed solely in terms of $N$:
\[
q_f(\Ff)\le Cq(\Ff)\max\{1,\log\log(N+2)\}.
\]
This consequence is stated explicitly in~\cite[Corollary~1.8]{pham_sharp_selectors}; the theorem and corollary numbers here refer to its September 2026 revision. The bounded-support result also improves the polynomial dependence on $t$ obtained independently by Dubroff, Kahn, and Park~\cite{dubroff2024note}.

Recently, Park~\cite[Theorem~1.2]{park} proved the dimension-free comparison
\[
q_f(\Ff)\le Kq(\Ff)
\max\{1,\log\log(1/q(\Ff))\}.
\]
Equivalently, every weakly $p$-small family is $(p/(K\max\{1,\log\log(1/p)\}))$-small up to changing $K$. Park's estimate implies the preceding dimension-dependent bound up to a universal constant, while its loss depends only on the density parameter. 

\begin{conjecture}[\cite{talagrand},Conjecture 6.3]\label{con:constant-bound}
There exist a universal constant $C$ such that for any increasing family $\Ff$,
\[
q_f(\Ff)\le Cq(\Ff).
\]
\end{conjecture}
This is the strongest form of this integral-to-fraction comparison. Equivalently, it conjectured that every weakly $p$-small family is $(p/C)$-small for a universal constant $C$. 

\begin{remark}
    On October 7, 2026, Conjecture \ref{con:constant-bound} is resolved by OpenAI in \cite{OAI:Integral-and-fractional-expectation-thresholds-are-equivalent-September-23-2026}. 
\end{remark}

Notably, this conjecture connects the study of random discrete structures with the rounding of fractional covers, which is another strong motivation to study this problem. Take any ground set $V$ and $p\in (0,1)$, let $\mu_p$ denote the product Bernoulli-$p$ measure on $2^V$. For every family $\Ff\subseteq 2^V$, define
\begin{equation}\label{def:F^m}
    \Ff^{(m)}=\{A\subseteq V:\forall S_1,\dots,S_m\in\Ff,A\not\subseteq S_1\cup\dots\cup S_m\}.
\end{equation}
A center conjecture of Talagrand is his discrete convexity conjecture \cite[Conjecture 7.1]{talagrand}, which states that there exists a integer $m$ such that for every ground set $V$, every $p\in(0,1)$ and every set family $\Ff$,
\[
\mu_p(\Ff)\ge 1-1/m\implies \Ff^{(m)} \text~{is}~p\text{-small}.
\]
Chen Li recently proved the fractional version in \cite[Corollary 3.1]{li}, which implies that
\begin{equation}\label{weakly-p-small}
    \mu_p(\Ff)>1/2\implies \Ff^{(2)} \text~{is weakly}~p/2\text{-small}.
\end{equation}
Thus, up to universal constants, the only remaining ingredient for integral convexity is Conjecture \ref{con:constant-bound}.

\subsection{Our Main Result and Some Definitions}

For \(x\geq1\), let \(\log_2^*x\) be the least integer \(k\geq0\)
such that applying \(\log_2\) repeatedly \(k\) times gives a number at
most \(1\). Each successive logarithm is taken only while the current
argument exceeds \(1\).

\begin{theorem}\label{thm:main}
Let \(V\) be a finite set, \(N=|V|\geq1\), and let
\(\mathcal{F}\subseteq2^V\) be a nontrivial increasing family. If
\(\mathcal{F}\) is weakly \(p\)-small for some \(0<p<1\), then it is
\[
    \frac{p}{64\log_2^*(N+2)}
\]
-small. Equivalently,
\[
    q_f(\mathcal{F})
    \leq64\log_2^*(N+2)\,q(\mathcal{F}).
\]
\end{theorem}

Using similiar methods, one can extend this result to the following stronger form: there exists a universal constant $C>0$ such that
    \begin{align}
    q_f(\Ff)\le Cq(\Ff)\log_2^*(l(\Ff)),\quad q_f(\Ff)\le Cq(\Ff)\bigl(1+\log_2^*(\operatorname{VC}(\min\Ff)+2)\bigr), \label{eq:ext}
    \end{align}
    Here $l(\Ff)$ is as defined above, and $\operatorname{VC}(\min\Ff)$ denotes the VC dimension of the minimal family. 
    The proofs for both refinements in \eqref{eq:ext} are given by AI. For completeness, a human-readable version of the proof for \eqref{eq:ext} is provided in Appendix~\ref{app:structural-refinements}.
    


Combining Theorem~\ref{thm:main} with \eqref{weakly-p-small} yields the following
corollary:

\begin{corollary}
If \(\mu_p(\mathcal{F})>1/2\), then
\[
    \mathcal{F}^{(2)}
    \text{ is }
    \frac{p}{128\log_2^*(N+2)}\text{-small}.
\]
\end{corollary}

We record notation used in the proofs. For \(U\subseteq V\), define
\[
    \langle U\rangle:=\{S\subseteq V:S\supseteq U\},
    \qquad
    \langle\mathcal{G}\rangle
    :=\{S\subseteq V:\text{some }G\in\mathcal{G}\text{ satisfies }G\subseteq S\}.
\]

\begin{definition}
For a nonnegative array
\(\lambda=(\lambda_S)_{\emptyset\ne S\subseteq V}\) and a nonempty
\(I\subseteq V\), set
\[
    P_\lambda(I):=\sum_{\emptyset\ne S\subseteq I}\lambda_S.
\]
We say that \(\lambda\) covers a family \(\mathcal{H}\) of nonempty
subsets of \(V\) if \(P_\lambda(I)\geq1\) for every \(I\in\mathcal{H}\).
\end{definition}

\begin{definition}
For \(r\in(0,1)\) and nonnegative array $\lambda$, define
\[
    A_\lambda(r):=\sum_{\emptyset\ne S\subseteq V}\lambda_Sr^{|S|},
    \qquad
    A_\lambda(r,U):=
    \sum_{\substack{S\subseteq V\\S\supseteq U}}\lambda_Sr^{|S|}
    \quad(\emptyset\ne U\subseteq V),
\]
and, for a family \(\mathcal{H}\) of nonempty subsets of \(V\), let
\[
    Z_r(\mathcal{H})
    :=\min_{\lambda\text{ covers }\mathcal{H}}A_\lambda(r).
\]
The minimum is attained because the problem has finitely many
variables and each objective coefficient is positive.
\end{definition}

Finally, if a weak-smallness certificate assigns weight
\(c=\lambda_\emptyset\) to the empty set, then \(c\leq1/2\).
Removing this weight and rescaling the others by \(1/(1-c)\)
preserves coverage and leaves the cost at most
\((1/2-c)/(1-c)\leq1/2\). We may therefore use certificates
supported only on nonempty sets.

\section{The Core-Extracting Argument}

\begin{lemma}[An efficient core]\label{lem:core}
Let $q,r,s\in(0,1)$ such that
\[
v:=r+q(1-r)\le s<1.
\]
Suppose a non-empty family $\Hh$ satisfies $Z_s(\Hh)\le B<1$. Choose
\[
\lambda\in\arg\min\{A_\lambda(r):\lambda~\text{covers}~\Hh\}.
\]
Then there is a nonempty $U$ with $A_\lambda(r,U)>0$ such that
\begin{equation}\label{eq:efficient}
q^{|U|}\le B\frac{A_\lambda(r,U)}{Z_r(\Hh)}.
\end{equation}
After buying $U$ and deleting all parents that contain $U$, the family $\Hh'=\Hh\setminus\langle U\rangle$ satisfies $Z_r(\Hh')\le Z_r(\Hh)-A_\lambda(r,U)$.
\end{lemma}

\begin{proof}
Let $\eta$ be a fractional cover of $\Hh$ such that $A_\eta(s)\leq B<1$. Let
\[
R=\max_{\varnothing\ne U\subseteq V}\frac{A_\lambda(r,U)}{q^{|U|}}.
\]
Since both $\lambda$ and $\eta$ covers $\Hh$, we have $P_\lambda(I)P_\eta(I)\ge 1$ for every $I\in \Hh$. Let $\gamma$ be another nonnegative array with $\gamma_W=\sum_{S\cup T=W}\lambda_S\eta_T$, then $P_\gamma(I)=P_\lambda(I)P_\eta(I)$, therefore it is consequently another fractional cover. Optimality of $\lambda$ gives
\begin{align}
Z_r(\Hh)
&=A_\lambda(r)\le A_\gamma(r)\notag\\
&=\sum_{S,T}\lambda_S\eta_Tr^{|S\cup T|}\notag\\
&=\sum_T\eta_Tr^{|T|}
\left[\sum_S\lambda_Sr^{|S|-|S\cap T|}\right]\notag\\
&=\sum_T\eta_Tr^{|T|}
\left[\sum_S\lambda_Sr^{|S|}\sum_{U\subseteq S\cap T}\left(\frac{1-r}{r}\right)^{|U|}\right]\notag\\
&=\sum_T\eta_Tr^{|T|}
\left[Z_r(\Hh)+\sum_{\varnothing\ne U\subseteq T}
\left(\frac{1-r}{r}\right)^{|U|}A_\lambda(r,U)\right]\label{eq:product}\\
&\le Z_r(\Hh)A_\eta(r)+R(A_\eta(v)-A_\eta(r)).\notag
\end{align}
The fourth identity follows by expanding
\[
r^{-|S\cap T|}
=\prod_{x\in S\cap T}\left(1+\frac{1-r}{r}\right)
=\sum_{U\subseteq S\cap T}\left(\frac{1-r}{r}\right)^{|U|}.
\]
For the last line, use $A_\lambda(r,U)\le Rq^{|U|}$ and
\[
r^{|T|}\left[\left(1+\frac{q(1-r)}r\right)^{|T|}-1\right]
=v^{|T|}-r^{|T|}.
\]
Since $A_\eta(x)$ is strictly increasing, we have $A_\eta(r)<A_\eta(v)\le A_\eta(s)\le B<1$. Therefore, by \eqref{eq:product} we have
\[
R\ge\frac{(1-A_\eta(r))Z_r(\Hh)}{A_\eta(v)-A_\eta(r)}\ge\frac{Z_r(\Hh)}{A_\eta(v)}\ge\frac{Z_r(\Hh)}{B}.
\]
Choose $U$ attaining $R$, we obtain $q^{|U|}\le B\frac{A_\lambda(r,U)}{Z_r(\Hh)}$. This proves \eqref{eq:efficient}.

Finally, notice that $\lambda':=(\lambda_S1_{S\not\supseteq U})_{\varnothing\ne S\subseteq V}$ is a fractional cover of $\Hh'$, we have $Z_r(\Hh')\le A_{\lambda'}(r)=Z_r(\Hh)-A_\lambda(r,U)$.
\end{proof}

\begin{lemma}[Logarithmic budget]\label{lem:budget}
Keep the larger-density cover $\eta$ and parameters of Lemma~\ref{lem:core} fixed. For any $0<\epsilon<Z_r(\Hh)$, there exists a family $\Gg$ such that 
\[
\sum_{S\in\Gg}q^{|S|}\le B\left(1+\log\frac{Z_r(\Hh)}{\epsilon}\right)
\]
and leave the residual family $\Hh':=\Hh\setminus\langle\Gg\rangle$ with $Z_r(\Hh')\le\epsilon$.
\end{lemma}

\begin{proof}
Let $\Hh_0=\Hh$. For every $n\ge 0$, while $Z_r(\Hh_n)>\varepsilon$, buy a core $U_n$ from Lemma~\ref{lem:core},adding $U_n$ to $\Gg$ and set $\Hh_{n+1}=\Hh_n\setminus\langle U_n\rangle$, and reoptimize. The fixed cover $\eta$ remains a fractional cover for every $\Hh_n$ with $s$-cost at most $B$. Each purchase deletes a parent, so the process is finite.

The cost of the $n$-th purchase is at most
\[
q^{|U_n|}\le B\frac{A_{\lambda_n}(r,U_n)}{Z_r(\Hh_n)}\le B\frac{Z_r(\Hh_n)-Z_r(\Hh_{n+1})}{Z_r(\Hh_n)}.
\]
For all steps before the last, $Z_r(\Hh_{n+1})>\varepsilon>0$, and
\[
\frac{Z_r(\Hh_n)-Z_r(\Hh_{n+1})}{Z_r(\Hh_n)}\le\log\frac {Z_r(\Hh_n)}{Z_r(\Hh_{n+1})}.
\]
These logarithms telescope. Thus all purchases before the last cost at most $B\log(Z_r(\Hh)/\varepsilon)$. The last purchase costs at most $B$, since $A_{\lambda_n}(r,U_n)\le Z_r(\Hh_n)$. This last-step bound also handles the possibility $Z_r(\Hh_{n+1})=0$, without taking a logarithm of zero. If the initial optimum is already at most $\varepsilon$, no purchase is needed.
\end{proof}

\begin{corollary}[One round]\label{cor:round}
Suppose $Z_s(\Hh)\le z<\beta$, with $z<1$, and $r=s-q>0$. At total integral $q$-cost at most $\beta$, one can leave a residual family $\Hh'$ satisfying
\begin{equation}\label{eq:round}
Z_r(\Hh')\le z\exp\left(1-\frac\beta z\right).
\end{equation}
\end{corollary}

\begin{proof}
Take an optimal cover at density $s$ as $\eta$. Then $B\le z$ and $Z_r(\Hh)\le z$. Also $r+q(1-r)\le r+q=s$. Set $\varepsilon=z\exp(1-\beta/z)<z$. If $Z_r(\Hh)\le\epsilon$, take $\Gg=\varnothing$ and $\Hh'=\Hh$. Otherwise, lemma~\ref{lem:budget} applies and bounds the total core cost by
\[
B\left(1+\log\frac z\varepsilon\right)
=B\frac\beta z\le\beta.
\]
\end{proof}

\section{Iteration and the Proof of Theorem \ref{thm:main}}
\begin{proof}[Proof of Theorem \ref{thm:main}]
    \textbf{Iteration of the one-round estimate.}
Take $\Hh_0=\min\Ff$; covering these parents covers all of $\Ff$. Since $\Ff$ is weakly $p$-small, $Z_p(\Hh_0)\le 1/2$. Then at
\[
s_0=\frac p{32},\qquad z_0=\frac1{64},
\]
its cost is at most $Z_{s_0}(\Hh_0)\le z_0$, since every witness is nonempty.

Let $d=\lstar(N+2)$ and $q=p/(64d)$. For $0\le j\le d$, set
\[
s_j=\frac p{32}-jq;
\qquad s_d=\frac p{64}\ge q.
\]
Define the numerical budget and cost sequences, for all $j\ge0$, by
\begin{equation}\label{eq:recurrence}
\beta_j=2^{-j-4},\qquad
z_{j+1}=z_j\exp\left(1-\frac{\beta_j}{z_j}\right).
\end{equation}
Starting with $\Hh_0$, apply Corollary~\ref{cor:round} at stage $j$ with old density $s=s_j$, new density $r=s-q=s_{j+1}$, and budget $\beta=\beta_j$. Reoptimization at the preceding stage provides the larger-density certificate required for the next stage. So once we have $z_j<\beta_j$, inductively this yields a residual family $\Hh_j$ with
\[
Z_{s_j}(\Hh_j)\le z_j.
\]
If a residual is empty, the construction has already finished.

To verify the condition $z_j<\beta_j$ and quantify the decay, write
\[
w_j=\frac{\beta_j}{2z_j}.
\]
Then $w_0=2$ and the exact recurrence is
\begin{equation}\label{eq:tower}
w_{j+1}=\frac{w_j}{2}\exp(2w_j-1)\ge\exp(w_j)
\qquad(w_j\ge2).
\end{equation}
The inequality follows from $(w_j/2)\exp(w_j-1)\ge1$. In particular, $w_j\ge2$ for every $j$, so $z_j\le\beta_j/4$ and every application of the corollary is valid.

For a fully explicit comparison with iterated logarithms, let
\[
t_0=1,\qquad t_{j+1}=2^{t_j}.
\]
The definition of $d$ implies $N+2\le t_d$. Since $w_0\ge t_0$ and $e^x\ge2^x$ for $x\ge0$, \eqref{eq:tower} implies $w_j\ge t_j$. Hence
\begin{equation}\label{eq:smallfinal}
w_d\ge N+2,\qquad
z_d=\frac{\beta_d}{2w_d}\le\frac1{32(N+2)}.
\end{equation}
The total cost of all purchased cores is less than
\[
\sum_{j\ge0}\beta_j=\frac18.
\]
\textbf{Final rounding, proved directly.}
If $\Hh_d$ is nonempty, let $\lambda$ be an optimal fractional cover at density $s_d$, with cost at most the bound in \eqref{eq:smallfinal}. Set
\[
a=(N+2)\log2.
\]
Independently select each nonempty witness $S$ with probability $1-e^{-a\lambda_S}$. For any parent $I\in\Hh_d$, its probability of being uncovered is
\[
\prod_{S\subseteq I}e^{-a\lambda_S}
=\exp\left(-a\sum_{S\subseteq I}\lambda_S\right)
\le e^{-a}.
\]
There are at most $2^N$ parents. By the union bound, the probability that some parent is uncovered is at most
\[
2^Ne^{-a}=\frac14.
\]
Because $q\le s_d$ and $1-e^{-x}\le x$ for $x\ge0$, the expected integral $q$-cost of the selected witnesses is at most
\[
\sum_Sq^{|S|}(1-e^{-a\lambda_S})
\le a\sum_S\lambda_Sq^{|S|}
\le az_d\le\frac{\log2}{32}.
\]
Markov's inequality shows that the probability this cost exceeds $1/4$ is at most $\log2/8$. As $1/4+\log2/8<1$, there is a realization covering every residual parent and having $q$-cost at most $1/4$.

Combine that realization with the previously purchased cores, removing duplicates if any. Every original parent is covered, and the total $q$-cost is less than $1/8+1/4=3/8$. This proves the covering assertion. Applying it to every weak-smallness parameter $p$ and taking the supremum of $p$ proves the threshold comparison.
\end{proof}

\section*{Statement on AI Usage and Human Effort} 


The proof in the main text was developed by GPT6 under the authors' guidance and was human-checked by the authors. More specifically, the authors proposed pursuing a deeper sunflower-type argument of Alweiss--Lovett--Wu--Zhang, which is shaped into the core extraction argument in its current form jointly by human and AI. Xuan supplied several candidate selection probabilities for the core extraction argument, and GPT6 refined the selection argument and finalized it in its current form. 

The Appendix contains a human-refactored presentation of the proofs for \eqref{eq:ext}; the original proofs in the Appendix were fully produced by GPT6. 




\clearpage
\appendix
\newtheorem{appendixtheorem}{\hspace{2em}Theorem}[section]
\newtheorem{appendixlemma}[appendixtheorem]{\hspace{2em}Lemma}
\numberwithin{equation}{section}

\section{Structural refinements}
\label{app:structural-refinements}

This appendix proves the two refinements stated in \eqref{eq:ext}, which are proved by AI and refactored into its current form by human.
Throughout, a parent means a member of $\min\Ff$, and all fractional
witnesses are nonempty, using the normalization in the introduction.
For an arbitrary family $\Hh$ of nonempty sets, we use $q(\Hh)$ and
$q_f(\Hh)$ for the thresholds of its upward closure.  Thus the same
covering definitions apply to parent families that are not increasing.
We set $Z_r(\varnothing)=0$ and write
\[
 T_0=1,\qquad T_{j+1}=2^{T_j}\quad(j\ge0).
\]
The definition of $\lstar$ gives $x\le T_{\lstar x}$ for $x\ge1$.
All logarithms without a specified base are natural.

\subsection{The largest minimal member}

\begin{appendixtheorem}[A parent-size comparison]
\label{thm:appendix-length}
Let $\Ff$ be a nontrivial increasing family on a finite ground set, and
put
\[
 l=l(\Ff)=\max\{2,\max_{I\in\min\Ff}|I|\}.
\]
If $\Ff$ is weakly $p$-small, then it is $p/(9\lstar l)$-small.
Consequently,
\[
 q_f(\Ff)\le9\lstar(l(\Ff))\,q(\Ff).
\]
\end{appendixtheorem}

\begin{proof}
Take a nonempty-witness certificate $\lambda^0$ with
$A_{\lambda^0}(p)\le1/2$.  Let $\mathcal A$ consist of all nonempty
unions of sets in its support, and define the restricted optimum
\[
 Z_r^{\mathcal A}(\Hh)
 =\min\{A_\lambda(r):\lambda\text{ covers }\Hh,
                  \ \lambda_S=0\text{ for }S\notin\mathcal A\}.
\]
The initial certificate is allowed.  Because $\mathcal A$ is closed
under unions, the product certificate in Lemma~\ref{lem:core} is also
allowed.  Deleting coefficients leaves the support in $\mathcal A$.
Thus Lemmas~\ref{lem:core} and~\ref{lem:budget}, and
Corollary~\ref{cor:round}, remain valid for these restricted optima.
The purchased cores need not themselves belong to $\mathcal A$.

Write $\Hh_0=\min\Ff$ and set
\[
 h=\max_{I\in\Hh_0}|\{S\in\mathcal A:S\subseteq I\}|,
 \qquad D=\lstar l,\qquad q=\frac p{9D}.
\]
Every parent contains an allowed witness, so $1\le h\le2^l-1$.
Since $l\ge2$, we have $D\ge1$.
For $0\le j\le D$, put
\[
 s_j=\frac{2p}{9}-jq,
 \qquad z_0=\frac19,\qquad
 \beta_j=\frac49\,9^{-j},\qquad
 z_{j+1}=z_j\exp(1-\beta_j/z_j).
\]
All allowed witnesses are nonempty, and hence
\[
 Z_{s_0}^{\mathcal A}(\Hh_0)
 \le A_{\lambda^0}(2p/9)
 \le\frac29A_{\lambda^0}(p)\le\frac19.
\]
Apply the restricted version of Corollary~\ref{cor:round} for $D$
rounds, unless the residual family becomes empty first.
To justify each round and estimate its outcome, define
$w_j=\beta_j/(2z_j)$.  Then
\[
 w_0=2,\qquad
 w_{j+1}=\frac{w_j}{9}\exp(2w_j-1)\ge2^{w_j}
 \quad(w_j\ge2).
\]
Indeed, the logarithm of the ratio of the last two expressions is
$\log(w/9)+(2-\log2)w-1$, which is positive at $w=2$ and increasing
for $w\ge2$.  In particular $z_j\le\beta_j/4<\beta_j$, so every round
is valid, and induction gives $w_j\ge T_{j+1}$.
Since $l\le T_D$, we obtain
\[
 w_D\ge T_{D+1}=2^{T_D}\ge2^l>h.
\]

The remaining family $\Hh_D$ therefore has a restricted certificate
$\lambda$ of cost at most $z_D$ at density $s_D=p/9=Dq$.
If this family is empty, no final purchase is needed.  Otherwise buy
every allowed witness with $\lambda_S\ge1/h$.  Each remaining parent
contains at most $h$ allowed witnesses whose coefficients sum to at
least one, so these purchases cover it.  Their total $q$-cost is at most
\[
 \sum_{\lambda_S\ge1/h}q^{|S|}
 \le hA_\lambda(q)
 \le\frac hD A_\lambda(s_D)
 \le\frac hD z_D
 \le\frac{\beta_D}{2D}
 \le\frac29\,9^{-D}.
\]
Here $A_\lambda(q)\le A_\lambda(s_D)/D$ uses nonempty witnesses.
The cores bought in the $D$ rounds cost at most
\[
 \sum_{j=0}^{D-1}\beta_j=\frac12(1-9^{-D}).
\]
The combined price is at most
$\frac12(1-9^{-D})+\frac29\,9^{-D}<1/2$, proving the covering claim.
Taking the supremum over the weak-smallness parameters $p$ proves the
threshold comparison.
\end{proof}

\subsection{VC dimension of the minimal family}

A set $D\subseteq V$ is shattered by $\Hh$ if
$\{I\cap D:I\in\Hh\}=2^D$.  The VC dimension
$\operatorname{VC}(\Hh)$ is the largest size of a shattered set.
We use two results of Shangguan~\cite[Corollaries~1.2 and~2.5]{shangguan2026vc}.
There are absolute constants $C_{\mathrm{VC}},C_{\mathrm{tr}}\ge1$
with the following properties.  First, every nonempty family $\Hh$ of
nonempty sets with $\operatorname{VC}(\Hh)\le d$, $d\ge1$, satisfies
\begin{equation}\label{eq:appendix-vc-input}
 q_f(\Hh)\le C_{\mathrm{VC}}q(\Hh)\log(d+1).
\end{equation}
Second, call a probability measure $\mu$ on $2^V$ $k$-spread if
$\Pr_{I\sim\mu}(T\subseteq I)\le k^{-|T|}$ for every nonempty $T$.
If $k\ge8$ and $\operatorname{VC}(\operatorname{supp}\mu)\le d$, then
\begin{equation}\label{eq:appendix-trace-input}
 \sum_{\substack{T\subseteq V\\|T|=j}}
       \Pr_{I\sim\mu}(T\subseteq I)^2
 \le(C_{\mathrm{tr}}d)^d(4/k)^j
 \qquad(j\ge1).
\end{equation}
These two results are the external inputs to the VC comparison; all
remaining steps are proved below.

\begin{appendixlemma}[The tower iteration at an arbitrary stopping time]
\label{lem:appendix-tower}
Suppose $Z_p(\Hh_0)\le1/2$, let $m\ge1$ be an integer, and choose
$q>0$ with $mq<p/32$.  One can buy cores of total $q$-cost less than
$1/8$ so that the residual family $\Hh$ satisfies
\[
 r=\frac p{32}-mq>0,\qquad Z_r(\Hh)\le\frac1{32T_m}.
\]
\end{appendixlemma}

\begin{proof}
Use the iteration from the proof of Theorem~\ref{thm:main}, with $m$
rounds and the present value of $q$.  Explicitly, start at
$s_0=p/32$ with cost at most $z_0=1/64$, and put
\[
 s_j=p/32-jq,\qquad \beta_j=2^{-j-4},\qquad
 z_{j+1}=z_j\exp(1-\beta_j/z_j).
\]
For $w_j=\beta_j/(2z_j)$, the recurrence \eqref{eq:tower} gives
$w_0=2$ and $w_{j+1}\ge2^{w_j}$, so $w_j\ge T_j$ and
$z_j<\beta_j$.  Corollary~\ref{cor:round} applies at each step.
Thus $Z_r(\Hh)\le z_m\le1/(32T_m)$, and the total cost is less than
$\sum_{j\ge0}\beta_j=1/8$.  An empty residual needs no further purchases.
\end{proof}

\begin{appendixlemma}[Counting the witnesses of an optimal certificate]
\label{lem:appendix-saturated}
Let $\Hh$ be a nonempty family of nonempty sets with
$\operatorname{VC}(\Hh)\le v$, where $v\ge1$ is an integer.
Suppose $0<r\le1/8$ and $z=Z_r(\Hh)<1$.  Put
\[
 B_v=(C_{\mathrm{tr}}(v+1))^{v+1}.
\]
For every optimal certificate $\lambda$ at density $r$,
\[
 |\{T:|T|=j,\ \lambda_T>0\}|\le B_v(4/r)^j
 \qquad(j\ge1).
\]
\end{appendixlemma}

\begin{proof}
Finite linear-programming duality gives an optimal nonnegative array
$(y_I)_{I\in\Hh}$ satisfying
\[
 \sum_{I\in\Hh}y_I=z,\qquad
 y(T):=\sum_{\substack{I\in\Hh\\I\supseteq T}}y_I\le r^{|T|}
 \quad(\varnothing\ne T\subseteq V).
\]
Complementary slackness implies
$\lambda_T>0\Rightarrow y(T)=r^{|T|}$.
For example, this follows directly by noting that both inequalities in
\[
 z=\sum_T\lambda_Tr^{|T|}
 \ge\sum_T\lambda_Ty(T)
 =\sum_Iy_IP_\lambda(I)\ge\sum_Iy_I=z
\]
must be equalities.

Give each $I\in\Hh$ probability mass $y_I$ and give $\varnothing$
the remaining mass $1-z$.  This defines a probability measure $\mu$
with $\Pr_{I\sim\mu}(T\subseteq I)=y(T)$ for nonempty $T$.
It is therefore $1/r$-spread.  Its support has VC dimension at most
$v+1$.  Indeed, if adjoining $\varnothing$ shatters a nonempty set
$D$, fix $x\in D$ and use the traces containing $x$ to see that the
original support shatters $D\setminus\{x\}$.
Since $1/r\ge8$, \eqref{eq:appendix-trace-input} gives
\[
 \sum_{|T|=j}y(T)^2\le B_v(4r)^j.
\]
Each $j$-set with $\lambda_T>0$ contributes $r^{2j}$ to this sum,
which proves the asserted count.
\end{proof}

\begin{samepage}
\begin{appendixtheorem}[An iterated-logarithm VC comparison]
\label{thm:appendix-vc}
There is an absolute constant $C$ such that every nontrivial increasing
family $\Ff$ satisfies
\[
 q_f(\Ff)\le
 C\bigl(1+\lstar(\operatorname{VC}(\min\Ff)+2)\bigr)q(\Ff).
\]
More precisely, writing $d=\operatorname{VC}(\min\Ff)$, weak
$p$-smallness implies $p/[C(1+\lstar(d+2))]$-smallness.
\end{appendixtheorem}
\end{samepage}

\begin{proof}
Let $p$ be a weak-smallness parameter and define
\[
 d=\operatorname{VC}(\min\Ff),\qquad v=\max\{1,d\},\qquad
 B=(C_{\mathrm{tr}}(v+1))^{v+1}.
\]
Choose
\[
 t=\lceil\log_8(16B)\rceil,\qquad
 K=\max\{2,2C_{\mathrm{VC}}\log(v+1)\},
\]
\[
 m=\lstar(K^t+2),\qquad q=\frac p{64(m+16)}.
\]
In particular $mq<p/32$.  Apply Lemma~\ref{lem:appendix-tower} to
$\min\Ff$.  After buying cores of total $q$-cost less than $1/8$, the
remaining actual parent family $\Hh\subseteq\min\Ff$ satisfies
\[
 r=\frac p{32}-mq=(m+32)q,\qquad
 z=Z_r(\Hh)\le\frac1{32T_m}
              \le\frac1{32(K^t+2)}.
\]
If $\Hh$ is empty, the construction is finished.  Otherwise choose an
optimal certificate $\lambda$ at density $r$.  We have $r\le p/32<1/8$
and $\operatorname{VC}(\Hh)\le v$, so
Lemma~\ref{lem:appendix-saturated} applies.

Buy every positive witness $S$ of $\lambda$ with $|S|>t$.
Since $4q/r\le1/8$, their total $q$-cost is at most
\[
 B\sum_{j>t}(4q/r)^j
 \le\frac{B8^{-t}}7\le\frac1{112}.
\]
Let $\mathcal J$ be the subfamily of $\Hh$ still uncovered by these
purchases.  Every member of $\mathcal J$ contains no positive long
witness, so the short certificate
$\lambda^{\mathrm{sh}}=(\lambda_S\mathbf1_{\{|S|\le t\}})_S$
still covers $\mathcal J$.  Moreover
$\operatorname{VC}(\mathcal J)\le v$, because we have only deleted
original parents.

Put $R=Kq$.  In the next estimate we evaluate $A_\lambda$ by its
defining finite polynomial also when its argument is at least one.
The short certificate satisfies
\[
 A_{\lambda^{\mathrm{sh}}}(R)
 \le\max\{1,R/r\}^{t}A_\lambda(r)
 \le K^t z\le\frac1{32}.
\]
If $\mathcal J$ is nonempty and $R\ge1$, its covering constraint at
any parent would force this polynomial to be at least one, a
contradiction.  Thus either $\mathcal J$ is empty, or $R<1$ and it is
weakly $R$-small.  In the latter case, \eqref{eq:appendix-vc-input}
implies
\[
 q(\mathcal J)\ge
 \frac{R}{C_{\mathrm{VC}}\log(v+1)}\ge2q.
\]
On a finite ground set the integral covering cost is the minimum of
finitely many continuous functions of the density, so the smallness
threshold is attained.  Hence $\mathcal J$ has an integral cover of
cost at most $1/2$ at density $2q$.  At density $q$ its nonempty
witnesses cost at most $1/4$.
Combining the three parts gives a cover of the original parents with
total $q$-cost less than
\[
 \frac18+\frac1{112}+\frac14=\frac{43}{112}<\frac12.
\]

It remains to bound the density loss.  We have
$t=O(v\log(v+1))$ and
$\log K=O(1+\log\log(v+2))$.  Consequently, for all sufficiently
large $v$,
\[
 \log_2(K^t+2)
 =O\bigl(v\log(v+1)(1+\log\log(v+2))\bigr)\le2^v.
\]
Thus $m\le\lstar v+2$ for those $v$; the bounded remaining values
are absorbed by an absolute additive constant.  Since $v=\max\{1,d\}$,
\[
 64(m+16)\le C(1+\lstar(d+2))
\]
for an absolute $C$.  Decreasing the density if necessary proves the
weak-smallness assertion.  Taking the supremum over $p$ proves the
threshold comparison.
\end{proof}

Both refinements concern the geometry of the original minimal family.
The first uses the number of allowed witnesses inside each parent;
the second uses VC estimates to pay for long optimal witnesses and
then rounds the remaining short certificate.  Neither argument requires
a bound on the number of vertices of the ground set.

\bibliographystyle{plain}
\bibliography{references}

\end{document}